\documentclass[hidelinks,onefignum,onetabnum]{siamart220329}

\usepackage{todonotes}
\usepackage{algorithmic}
\usepackage{algorithm}
\usepackage{amsmath}
\usepackage{amssymb}
\usepackage{graphicx}
\usepackage{listings}

\usepackage{wrapfig}
\usepackage{tikz}

\newtheorem{problem}{Problem}

\newcommand{\R}{{\mathbb R}}

\newcommand{\ass}{\,\mbox{:=}\,}

\newcommand{\fma}{{\footnotesize \,\text{fma}}}

\usepackage{graphicx}

\begin{document}

\headers{A Matrix-Free Newton Method}{U. Naumann}

\title{A Matrix-Free Newton Method}

\author{Uwe Naumann\thanks{Software and Tools for Computational Engineering,
               RWTH Aachen University,
               Aachen, Germany,
	      \email{naumann@stce.rwth-aachen.de}}}

\maketitle

\begin{abstract}
A modification of Newton's method for solving systems of $n$ nonlinear 
equations is presented. The new matrix-free method relies on a given 
decomposition of the invertible Jacobian of the residual into invertible 
sparse local Jacobians according to the chain rule of differentiation.
It is motivated in the context of local Jacobians with bandwidth $2m+1$ for 
$m\ll n$. A reduction of the computational cost by $\mathcal{O}(\frac{n}{m})$ 
can be observed. Supporting run time measurements are presented for the 
tridiagonal case showing a reduction of the computational cost by $\mathcal{O}(n).$ 

Generalization yields the combinatorial {\sc Matrix-Free Newton Step} problem.
We prove NP-completeness and we present algorithmic components for
building methods for the approximate solution. Inspired by adjoint Algorithmic
Differentiation, the new method shares several challenges for the latter 
including the {\sc DAG Reversal} problem. Further
challenges are due to combinatorial problems in sparse linear algebra such as
{\sc Bandwidth} or {\sc Directed Elimination Ordering}.
\end{abstract}



\section{Introduction} \label{sec:0}

We revisit Newton's method \cite{Newton1760Pnp,Whiteside1976TMP} for computing roots 
$x \in \R^n$ of differentiable multivariate vector functions 
$F : \R^n \rightarrow \R^n : y=F(x)$ 
with invertible dense Jacobians
$$
F'(x) \equiv \left (\frac{d y_j}{d x_i} \right)_{i,j=1,\ldots,n} \in \R^{n \times n} \; ,
$$
and where $x=(x_i)$ and $y=(y_i)$ for $i=1,\ldots,n.$ Approximate solutions
for which $F(x) \approx 0$ are computed iteratively as
\begin{equation} \label{eqn:newton}
\begin{split}
\Delta x&\ass -F'(x)^{-1} \cdot F(x) \\
x&\ass x+\Delta x
\end{split}
\end{equation}
for given starting points $x \in \R^n.$ 
Convergence after $p\geq 0$ 
iterations 
is typically defined as the 
norm of the residual $F(x)$ falling below some given threshold $0<\epsilon \ll 1.$ 
See, for example, \cite{Deuflhard2004NMf,Kelley2003SNE} for further information on Newton's method.

We use $=$ to denote mathematical equality, $\equiv$ in the sense of ``is 
defined as'' and $\ass$ to represent assignment according to imperative 
programming. Approximate equality if denoted as $\approx$.
We distinguish between partial ($\partial$) and total ($d$) derivatives.
Multiplication is denoted by a dot. The dot may also be omitted in favor of a 
more compact notation.

The residual $F$ is usually given as a differentiable program written in 
some high-level programming language.
Equation~(\ref{eqn:newton}) implies that the two main ingredients of Newton's 
method are
the evaluation of the Jacobian $F'(x)$ and
the solution of the linear Newton system $F'(x) \cdot \Delta x =-F(x).$
Algorithmic Differentiation (AD) \cite{Griewank2008EDP} of $F$ 
yields $F'(x)$ with machine accuracy at a
computational cost of ${\mathcal O}(n) \cdot \text{Cost}(F).$ 
Subsequent direct solution of the linear Newton system
yields the Newton steps $\Delta x \in \R^n$ at the
computational cost of ${\mathcal O}(n^3).$ 
The overall computational cost can be dominated by either of the two
parts depending on the ratio $\frac{\text{Cost}(F)}{{\cal O}(n^2)}.$  

Exploitation of special structure of $F$ can yield a
significant reduction in computational cost. 
Consequently, this article proposes a {\em matrix-free Newton method} motivated
by savings obtained for relevant scenarios. We present the 
fundamental idea behind the method and we discuss further generalization.
In Section~\ref{sec:1} we recall essential 
fundamentals of AD and we draw conclusions for the computation of Newton 
steps. The potential for reduction in computational cost is illustrated
with the help of the practically relevant special case of banded local Jacobians in Section~\ref{sec:2}. Further generalization and formalization is the subject of Section~\ref{sec:3}. Conclusions are drawn in Section~\ref{sec:4}.

\section{A Lesson from Adjoint AD} \label{sec:1}

First-order AD comes in two fundamental flavors. Tangent AD yields
\begin{equation} \label{tad}
	\dot{y}=\dot{F}(x,\dot{x}) \equiv F'(x) \cdot \dot{x}
\end{equation}
and, hence, the (dense) Jacobian with machine accuracy at ${\cal O}(n) \cdot \text{Cost}(\dot{F})$ 
by letting $\dot{x}$ range over the Cartesian basis of $\R^n.$
Adjoint AD computes 
\begin{equation} \label{aad}
	\bar{x}=\bar{F}(x,\bar{y}) \equiv F'(x)^T \cdot \bar{y} \; ,
\end{equation}
yielding the 
same Jacobian up to machine accuracy at ${\cal O}(n) \cdot \text{Cost}(\bar{F})$ 
by letting $\bar{y}$ range over the Cartesian basis of $\R^n.$ 
Both tangent and adjoint AD are matrix-free methods in the sense that
the Jacobian is not required explicitly in order to evaluate Equations~(\ref{tad})
or (\ref{aad}). Both methods can be used to accumulate $F'.$
Their costs differ according to the ratio
$$
{\cal R} \equiv \frac{\text{Cost}(\bar{F})}{\text{Cost}(\dot{F})} \; ,
$$
where, typically, ${\cal R} \ge 1.$ Hence, tangent AD is usually the method
of choice for computing Jacobians of the residual in the context of Newton's 
method. Sparsity of the Jacobian of the residual may change the picture \cite{Gebremedhin2005WCI}.

Second-order tangent and adjoint AD follow naturally; 
see \cite{Griewank2008EDP}. Hessians required by Newton's method for convex
optimization of objectives $f : \R^n \rightarrow \R$ can be computed at
${\cal O}(n) \cdot \text{Cost}(\dot{\bar{f}})$ using second-order
adjoint versions $\dot{\bar{f}}$ of $f.$ Corresponding matrix-free 
Newton-Krylov methods \cite{Knoll2004JfN} can be derived. They are based on the observation that the 
Hessian-vector products required by Krylov-subspace methods (e.g. Conjugate 
Gradients \cite{Hestenes1952MoC}) can be computed by a second-order 
adjoint $\dot{\bar{f}}$ without
prior accumulation of the Hessian. The approach to be proposed in this article
is different. It exploits special structure and local sparsity of $F$ instead
of global propagation of tangents or adjoints.

\subsection{Terminology}

AD requires the given program $F$ to be differentiable. The notation from \cite{Griewank2008EDP} is modified only slightly.
\begin{definition}[Differentiable Program] \label{def:diffprog}
	A {\em differentiable program} $$F : \R^n \rightarrow \R^m : \; y\ass F(x)$$ 
	decomposes into a {\em single assignment code} 
	\begin{equation} \label{eqn:sac}
	v_j \ass \varphi(v_i)_{i \in P_j} \quad \text{for}~j=1,\ldots,p+m \; ,
	\end{equation}
where $v_{i-n}=x_i$ for $i=1,\ldots,n$ and $y_k=v_{p+k}$ for $k=1,\ldots,m$
	and with differentiable {\em elemental functions} $\varphi_j,$ $j=1,\ldots,p+m,$ 
	featuring {\em elemental partial derivatives} 
$$
	\partial_{j,i} \equiv \frac{\partial \varphi_j}{\partial v_i} \quad
	\text{for}~i \in P_j \; .
$$
	The set of arguments of $\varphi_j$ (direct predecessors of $v_j$) is denoted by $P_j.$
\end{definition}

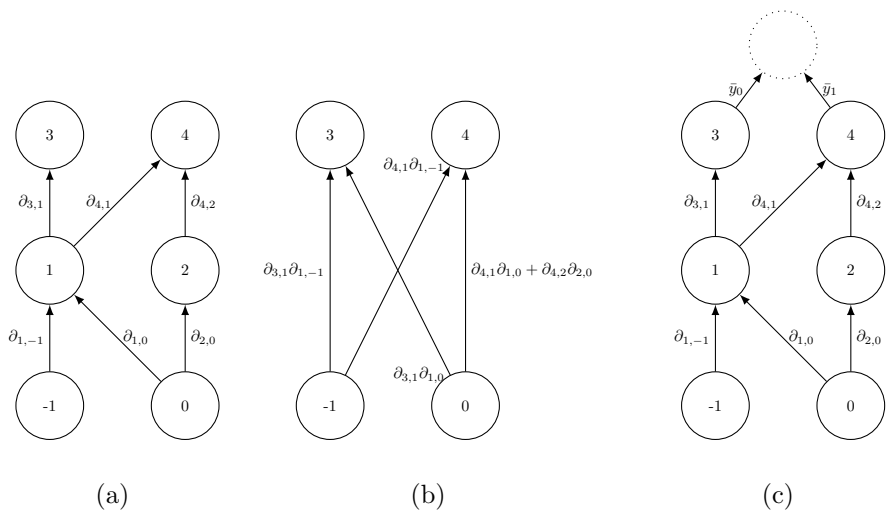
\begin{figure}
	    \begin{tabular}{ccc}
                \begin{minipage}[b]{.23\textwidth}
\begin{center}
\begin{tikzpicture}[scale=.6, transform shape]
            \begin{pgfscope}
                \tikzstyle{every node}=[draw,circle,minimum size=1.5cm]
                    \node (00) at (0,0) {-1};
                    \node (01) at (3,0) {0};
                    \node (10) at (0,3) {1};
                    \node (11) at (3,3) {2};
                    \node (20) at (0,6) {3};
                    \node (21) at (3,6) {4};
            \end{pgfscope}
            \begin{scope}[-latex]
                \draw (00) -- (10) node[midway,left] {$\partial_{1,-1}$};
                \draw (01) -- (11) node[midway,right] {$\partial_{2,0}$};
                \draw (01) -- (10) node[midway,right] {$\partial_{1,0}$};
                \draw (10) -- (20) node[midway,left] {$\partial_{3,1}$};
                \draw (11) -- (21) node[midway,right] {$\partial_{4,2}$};
                \draw (10) -- (21) node[midway,left] {$\partial_{4,1}$};
            \end{scope}
        \end{tikzpicture}
\end{center}
                \end{minipage}
		    &
            \begin{minipage}[b]{.35\textwidth}
\begin{center}
\begin{tikzpicture}[scale=.6, transform shape]
            \begin{pgfscope}
                \tikzstyle{every node}=[draw,circle,minimum size=1.5cm]
                    \node (00) at (0,0) {-1};
                    \node (01) at (3,0) {0};
                    \node (20) at (0,6) {3};
                    \node (21) at (3,6) {4};
            \end{pgfscope}
            \begin{scope}[-latex]
                \draw (00) -- (20) node[midway,left] {$\partial_{3,1} \partial_{1,-1}$};
                \draw (00) -- (21) node[at end,left] {$\partial_{4,1} \partial_{1,-1}$};
                \draw (01) -- (20) node[at start,left] {$\partial_{3,1} \partial_{1,0}$};
                \draw (01) -- (21) node[midway,right] {$\partial_{4,1} \partial_{1,0}+\partial_{4,2} \partial_{2,0}$};
            \end{scope}
        \end{tikzpicture}
\end{center}
                \end{minipage} 
                &
		    \begin{minipage}[b]{.3\textwidth}
        \centering
\begin{tikzpicture}[scale=.6, transform shape]
            \begin{pgfscope}
                \tikzstyle{every node}=[draw,circle,minimum size=1.5cm]
                    \node (00) at (0,0) {-1};
                    \node (01) at (3,0) {0};
                    \node (10) at (0,3) {1};
                    \node (11) at (3,3) {2};
                    \node (20) at (0,6) {3};
                    \node (21) at (3,6) {4};
                    \node[dotted] (3) at (1.5,8) {};
            \end{pgfscope}
            \begin{scope}[-latex]
                \draw (00) -- (10) node[midway,left] {$\partial_{1,-1}$};
                \draw (01) -- (11) node[midway,right] {$\partial_{2,0}$};
                \draw (01) -- (10) node[midway,right] {$\partial_{1,0}$};
                \draw (10) -- (20) node[midway,left] {$\partial_{3,1}$};
                \draw (11) -- (21) node[midway,right] {$\partial_{4,2}$};
                \draw (10) -- (21) node[midway,left] {$\partial_{4,1}$};
                    \draw (20) -- (3) node[midway,left] {$\bar{y}_0$};
                    \draw (21) -- (3) node[midway,right] {$\bar{y}_1$};
            \end{scope}
        \end{tikzpicture}
\end{minipage}
		    \\
		    \\
		    (a) & (b) & (c)
\end{tabular}
	\caption{Sample DAGs} \label{fig:1}
\end{figure}

\begin{definition}[Direct Acyclic Graph] \label{def:dag}
	Equation~(\ref{eqn:sac}) induces a {\em directed acyclic graph (DAG)} 
$G=(V,E),$ $V=(X,Z,Y),$ $X \cap Z \cap Y=\varnothing,$\footnote{The 
simplest ``no-op'' program \lstinline{x=x} distinguishes \lstinline{x} as 
an input ($\in X$) from \lstinline{x} as an output $(\in Y)$ while 
$Z=\varnothing$ in this case. Corresponding larger scenarios follow naturally.} $E \subseteq V \times V$ such that $X=\{1-n,\ldots,0\},$ $Z=\{1,\ldots,p\},$
$Y=p+1,\ldots,p+m$ and $(i,j) \in E \Leftrightarrow i \in P_j.$ 
All edges $(i,j)$ are labelled with $\partial_{j,i}.$ 
\end{definition}
Paths connecting vertices $i$ and $j$ are denoted as $(i,\ldots,j).$
\begin{definition}[Layered DAG and Local Jacobians] \label{def:ldag}
A DAG of a matrix chain product 
$$
F'_q \cdot \ldots \cdot F'_1 = F' \in \R^{m \times n}
$$
	of {\em local Jacobians} $F'_i,$ $i=1,\ldots,q$ is called {\em layered DAG}.
\end{definition}
The DAG in Figure~\ref{fig:1}(a) is layered.
\begin{definition}[Uniformly Layered DAG] \label{def:uldag}
A layered DAG of a matrix chain product
$$
F'_q \cdot \ldots \cdot F'_1 = F' \in \R^{n \times n}
$$
	is called {\em uniformly layered DAG} if 
$F'_i \in \R^{n \times n}$ for $i=1,\ldots,q.$ 
\end{definition}
The DAG in Figure~\ref{fig:1}(a) is uniformly layered.
\begin{definition}[Invertible DAG] \label{def:idag}
We refer to a DAG of a differentiable program with invertible Jacobian 
$F' \in \R^{n \times n}$ as an {\em invertible DAG}. 
The {\em corresponding Newton step} is 
defined as $(F')^{-1} \cdot y$ for given $y \in \R^n.$
\end{definition}

Adjoint AD can be implemented by using operator and function overloading in
suitable programming languages such as C++ \cite{Griewank1996AAC,Hogan2014FRM,Sagebaum2019HPD}. 
The given implementation of $F$
as a differentiable program is run in overloaded arithmetic to augment the
computation of the function value with the recording of the corresponding DAG
(also referred to as {\em tape}). 
For example, for \lstinline{q=2} the differentiable C++ program
\begin{lstlisting}
for (int i=0;i<q;i++) {
  x[i%2]*=x[(i+1)%2];
  x[(i+1)%2]=sin(x[(i+1)%2]);
}
\end{lstlisting}
yields the DAG in Figure~\ref{fig:1}(a), where, $\partial_{1,-1}=\text{x[1]},$ $\partial_{3,1}=\cos(\text{x[0]}),$ and so forth. The Jacobian 
$$
F'=
\begin{pmatrix}
	\partial_{3,1} \partial_{1,-1} &  \partial_{3,1} \partial_{1,0} \\
	\partial_{4,1} \partial_{1,-1} & \partial_{4,1} \partial_{1,0}+\partial_{4,2} \partial_{2,0}
\end{pmatrix} = 
\begin{pmatrix}
	\partial_{3,1} & 0 \\
	\partial_{4,1} & \partial_{4,2}
\end{pmatrix} 
\begin{pmatrix}
	\partial_{1,-1} & \partial_{1,0} \\
	0 & \partial_{2,0}
\end{pmatrix} 
$$
can be represented as the bipartite DAG in Figure~\ref{fig:1}(b). 
Both options for evaluating the adjoint 
$$
\left (
\begin{pmatrix}
	\partial_{3,1} & 0 \\
	\partial_{4,1} & \partial_{4,2}
\end{pmatrix} 
\begin{pmatrix}
	\partial_{1,-1} & \partial_{1,0} \\
	0 & \partial_{2,0}
\end{pmatrix} 
\right )^T 
\begin{pmatrix}
	\bar{y}_0 \\
	\bar{y}_1 \\
\end{pmatrix} =
\begin{pmatrix}
	\partial_{1,-1} & 0 \\
	\partial_{1,0} & \partial_{2,0}
\end{pmatrix}
\left (
\begin{pmatrix}
	\partial_{3,1} & \partial_{4,1} \\
	0 & \partial_{4,2}
\end{pmatrix} 
\begin{pmatrix}
	\bar{y}_0 \\
	\bar{y}_1 \\
\end{pmatrix} 
\right )
$$
(with $\bar{y}=(0~1)^T$ and $\bar{y}=(1~0)^T$ giving the rows of the Jacobian) yield the same result. Note
that the evaluation of the expression on the left-hand side requires 12\fma\footnote{fused multiply-add floating-point operations} (${\cal O}(q   n^3)$; here $q=n=2$) while
the right-hand-side expression take only 8\fma~(${\cal O}(q   n^2)$).
The same effect due to associativity of the chain rule of differentiation
(equivalently, of matrix multiplication)\footnote{assuming infinite precision arithmetic} is exploited by {\em backpropagation} in the context of training 
of artificial neural networks \cite{GoodBengCour16}. The adjoint can be 
visualized as the DAG in Figure~\ref{fig:1}(c).

A discussion of the numerous aspects of AD and of its implementation are beyond 
the scope of this article. The interested reader is referred to 
\cite{Griewank2008EDP} for further information on the subject. 
Moreover, the AD community's
Web portal {\tt www.autodiff.org} contains a comprehensive bibliography as
well as links to various AD research and tool development efforts.

\subsection{Lesson}

\begin{wrapfigure}[15]{l}{4cm}
	\centering
\begin{tikzpicture}[scale=.6, transform shape]
            \begin{pgfscope}
                \tikzstyle{every node}=[draw,circle,minimum size=1.2cm]
                    \node (00) at (0,0) {};
                    \node (01) at (2,0) {};
                    \node (10) at (0,2) {};
                    \node (11) at (2,2) {};
                    \node (20) at (0,4) {};
                    \node (21) at (2,4) {};
                    \node[dotted] (3) at (1,6) {};
            \end{pgfscope}
            \begin{scope}[-latex]
                \draw (10) -- (00) node[midway,left] {$a$};
                \draw (11) -- (01) node[midway,right] {$b$};
                \draw (10) -- (01) node[midway,left] {$c$};
                \draw (20) -- (10) node[midway,left] {$d$};
                \draw (21) -- (11) node[midway,right] {$e$};
                \draw (21) -- (10) node[midway,left] {$f$};
                    \draw (3) -- (20) node[midway,left] {$y_0$};
                    \draw (3) -- (21) node[midway,right] {$y_1$};
            \end{scope}
        \end{tikzpicture}
	\caption{Inverse DAG} \label{fig:2}
\end{wrapfigure}
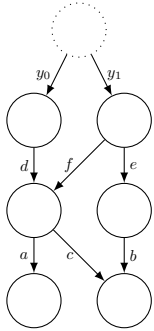
How can Newton's method benefit from lessons learned in adjoint AD? Working 
towards an answer to this question,
the given implementation of $F$ as a differentiable program is assumed to be composed of elemental 
functions 
$F_i : \R^n \rightarrow \R^n : x^i=F_i(x^{i-1})$ for $i=1,\ldots,q$
as
\begin{equation} \label{eqn:1}
y=x^q=F(x^0)=(F_q  \circ \ldots \circ F_1)(x^0)
\end{equation} 
where $x^0=x.$ 
Application of the chain rule of differentiation to 
Equation~(\ref{eqn:1}) yields 
\begin{equation} \label{eqn:2}
F'(x^0)=F'_q(x^{q-1}) \cdot \ldots \cdot F'_1(x^0)
\end{equation} implying 
\begin{equation} \label{eqn:3}
	\begin{split}
		\Delta x&=
		-F'(x)^{-1} \cdot y \\
		&= -(F'_q(x^{q-1}) \cdot \ldots \cdot F'_1(x^0))^{-1} \cdot y \\
		&=-F'_1(x^0)^{-1} \cdot \ldots \cdot F'_q(x^{q-1})^{-1} \cdot y\; .
	\end{split}
\end{equation}
$F'(x)$ needs to be invertible. Invertibility of all $F'_i$ follows immediately.
Their (w.l.o.g.) $LU$ factorization \cite{Duff1986DMf} yields 
\begin{equation} \label{eqn:3.1}
F_i \equiv F_i(x^{i-1})=L_i(x^{i-1}) \cdot U_i(x^{i-1}) \equiv 
L_i \cdot U_i
\end{equation}
with lower triangular $L_i$ and upper triangular $U_i.$ 
In-place factorization yields $L_i-I_n+U_i$ 
at a computational cost of ${\mathcal O}(n^3).$ 
The matrix $I_n \in \R^{n \times n}$ denotes the identity in $\R^n.$ 
From Equation~(\ref{eqn:3}) it follows that
\begin{equation} \label{eqn:4}
\begin{split}
\Delta x&=-U_1^{-1} \cdot L_1^{-1} \cdot \ldots \cdot U_q^{-1} \cdot L_q^{-1} \cdot y \\
&=-U_1^{-1} \cdot \left ( L_1^{-1} \cdot \ldots \cdot\left (  U_q^{-1} \cdot \left ( L_q^{-1} \cdot y \right )\right ) \ldots \right ) 
\end{split}
\end{equation}
yielding $2  q$ linear systems to be solved as efficiently as possible in 
order to undercut the computational cost of the standard Newton method.
The new method is matrix-free in the sense that $F'$ is not accumulated explicitly.
The analogy with adjoint AD is illustrated by 
Figure~\ref{fig:2}, where the reversed edges indicate products of vectors 
with inverse local Jacobians. For the example in Figure~\ref{fig:1}, the standard {\em ``accumulate first, then factorize''} approach yields the Newton step
$$
        \underset{\Sigma=12}{\underbrace{
        \underset{{\cal O}(n^3)~\hat{=}~3+3+1=7}{\underbrace{
        \underset{{\cal O}(q   n^3)~\hat{=}~5}{\underbrace{
        \left (
\begin{pmatrix}
	\partial_{3,1} & 0 \\
	\partial_{4,1} & \partial_{4,2}
\end{pmatrix} 
\begin{pmatrix}
	\partial_{1,-1} & \partial_{1,0} \\
	0 & \partial_{2,0}
\end{pmatrix} 
        \right )^{-1}
        }}
        \begin{pmatrix}
                y_0 \\
                 y_1
        \end{pmatrix} }}}}
        $$
at the expense of $12\fma.$ The alternative {\em ``factorize first, then accumulate''} method   
$$
        \underset{{O}(q   n^2)~\hat{=}~{6}}{\underbrace{
\begin{pmatrix}
	\partial_{1,-1} & \partial_{1,0} \\
	0 & \partial_{2,0}
        \end{pmatrix}^{-1} 
        \underset{{O}(n^2)~\hat{=}~3}{\underbrace{
        \left (
\begin{pmatrix}
	\partial_{3,1} & 0 \\
	\partial_{4,1} & \partial_{4,2}
        \end{pmatrix}^{-1} 
        \begin{pmatrix}
                y_0 \\
                 y_1 
        \end{pmatrix} \right )}}}}
       $$
performs the same task using only $6\fma.$ Admittedly, this example is extreme
in the sense that all elemental Jacobians are already triangular. Products of their
inverses with
a vector can hence be computed very efficiently by simple forward or backward 
substitution. More realistically, the individual factors need to be transformed
into triangular form first. Ideally, the exploitation of sparsity of the elemental 
Jacobians is expected to keep the corresponding additional effort 
low. However, the usual challenges faced in the context of direct methods for
sparse linear algebra \cite{Duff1986DMf} need to be addressed. A representative special case is 
discussed in the next section. It illustrates the potential of the new method 
and it serves as motivation for further generalization in Section~\ref{sec:3} 
with the aim to enable applicability to a wider range of practically relevant 
problems.

\section{Banded Elemental Jacobians} \label{sec:2}

Let all $F'_i$ have (maximum) bandwidth $2m+1$ ($m$ off-diagonals)
with $m \ll n.$ Their in-place $LU$ factorization as in Equation~(\ref{eqn:3.1})
yields $L_i-I_n+U_i$ with the same bandwidth at a computational cost of 
${\mathcal O}(m^2  n).$
The resulting $2   q$ triangular linear systems can be solved efficiently 
with a cost of ${\mathcal O}(m   n)$ by simple substitution, respectively.
The total computational cost of the matrix-free Newton method can hence be 
estimated as ${\mathcal O}(q   m^2   n).$ 

For $F'$ to become dense we require $q \ge \underline{q} \equiv (n-m-1)/m.$ 
For $q\gg \underline{q}$ the cost of computing $F'$ can be estimated as
${\mathcal O}(q   m   n^2).$ A matrix chain product of length $q$ 
needs to be evaluated, where one factor has bandwidth $2m+1$ and the other 
factor becomes dense for $q\ge \underline{q}.$
Superiority of the matrix-free Newton method
follows immediately for $m \ll n$ as
${\mathcal O}(q   m^2   n) < {\mathcal O}(q   m   n^2).$ 

In most 
real-world scenarios the accumulation of $F'$ is likely to dominate the overall
computational cost yielding a reduction 
of the computational cost of the matrix-free Newton method over the cost of
the standard Newton method by $\mathcal{O}(\frac{n}{m}).$ 
Nevertheless, let the computational cost of a Newton step be dominated by 
the solution of the linear Newton system. 
Superiority of the matrix-free Newton method requires ${\cal O}(q  m^2  n)<{\cal O}(n^3),$ which will only be violated for $q>n$ assuming
$n={\cal O}(m^2).$ 

\begin{figure}
	\includegraphics[width=\textwidth]{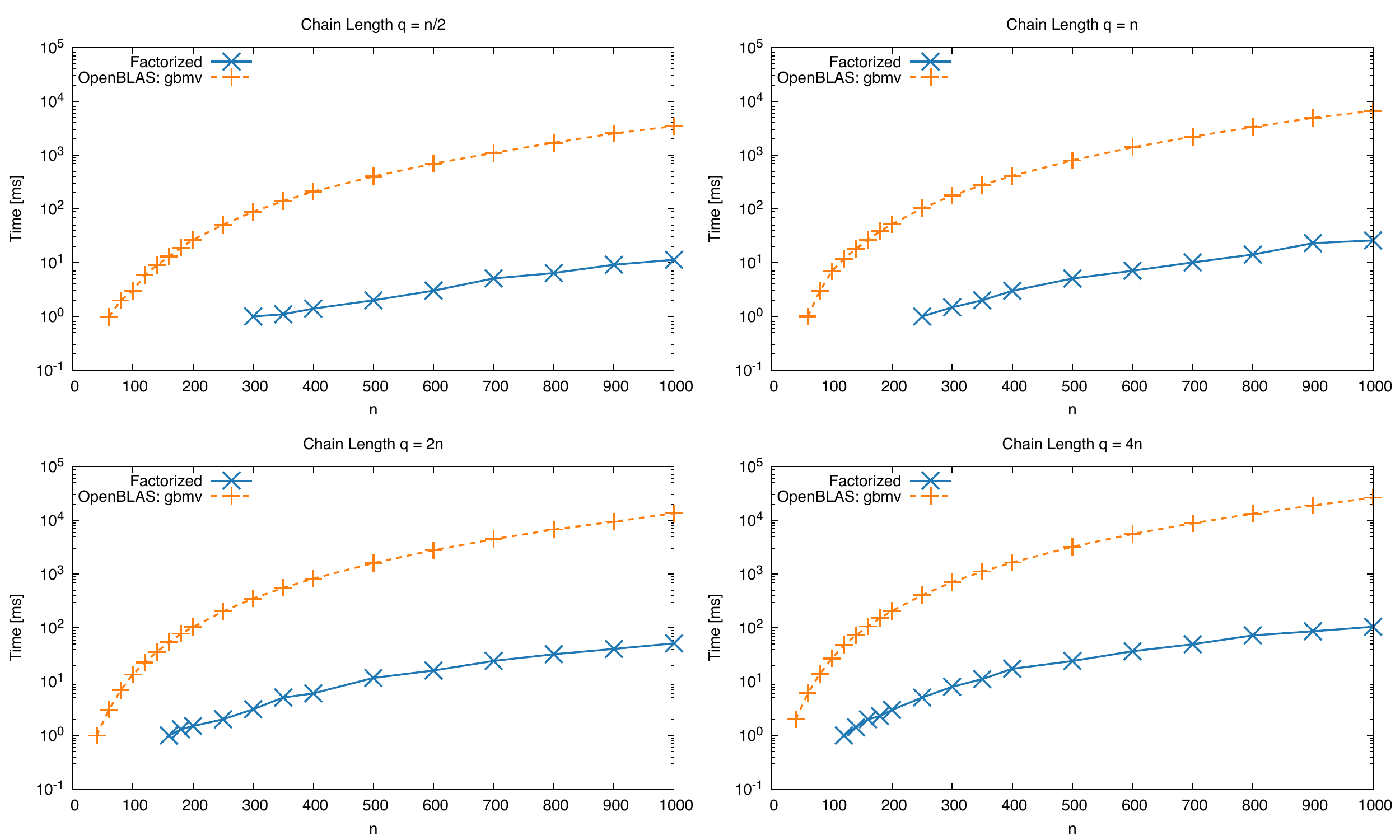}
	\caption{Tridiagonal Elemental Jacobians} \label{fig:3}
\end{figure}

Let all $F'_i$ be tridiagonal. For example, this scenario occurs in the context of implicit Euler integration for one-dimensional diffusion with spacial
discretization performed by central finite differences.
In-place $LU$ factorization using the Thomas algorithm \cite{Thomas1949EPi} 
yields the tridiagonal
$L_i-I_n-U_i$ at a computational cost of ${\mathcal O}(n).$ 
The $2   q$ linear systems in Equation~(\ref{eqn:4}) can be solved 
with a cost of the same order, respectively.
The total computational cost of the matrix-free Newton method adds up to 
${\mathcal O}(q   n) < \max({\mathcal O}(q   n^2), {\mathcal O}(n^3)).$ 
The speedup of $\mathcal{O}(\frac{n}{m})=\mathcal{O}(n)$ is illustrated in 
Figure~\ref{fig:3} showing the results of various experiments for a single
Newton step and 
$q/n=0.5,1,2,4,$ $0<n\le 10^3.$ The 
{\em ``factorize first, then accumulate''} approach implements the new
matrix-free Newton method. It yields the bottom line in all four sub-figures of Figure~\ref{fig:3}.
The {\em ``accumulate first, then factorize''} approach uses OpenBLAS' 
{\tt gbmv} method ({\tt www.openblas.net}) for multiplying banded matrices 
with dense vectors. It exceeds the computational cost of the new method by 
roughly a factor of $n.$ Our search
for a dedicated method for multiplying banded matrices turned out unsuccessful.
We would expect such an algorithm to be more efficient for $q \ll n.$ Our 
basic reference implementation did not outperform {\tt gbmv} though. Little
effort went into its optimization as it would not be our method of choice for
the more common scenario of $q \gg n.$

\section{Towards Generalization} \label{sec:3}

DAGs recorded by adjoint AD are typically neither layered nor 
uniform. According to Definition~\ref{def:ldag} the term ``layered'' 
refers to a decomposition of the DAG 
into a sequence of bipartite sub-DAGs representing the elemental Jacobians 
$F'_i.$ For a DAG to be uniformly layered all layers 
must consist of the same number ($n$) of vertices as formally stated in
Definition~\ref{def:uldag}. 

In the following we present ideas on how to make invertible DAGs 
uniformly layered. We formalize the corresponding methods and we 
prove their numerical correctness, respectively. All proofs rely on the 
following formulation of the chain rule of differentiation.
\begin{lemma}[Chain Rule of Differentiation on DAG] \label{lem:cr}
Let $G$ be the DAG of a differentiable program $F : \R^n \rightarrow \R^m$ 
as in Definition~\ref{def:dag} with Jacobian 
	$$F' \equiv \left ( \frac{d y_k}{d x_i} \right )_{i=1,\ldots,n}^{k=1,\ldots,m} = \left ( \frac{d v_{p+k}}{d v_{i-n}} \right )_{i=1,\ldots,n}^{k=1,\ldots,m}
	\in \R^m \times \R^n \; .$$
	Then
	\begin{equation} \label{eqn:cr}
	\frac{d v_{p+k}}{d v_{i-n}} = \sum_{(i-n,\ldots,p+k)} \prod_{(s,t) \in (i-n,\ldots,p+k)} \partial_{t,s}
	\end{equation}
	for $i=1,\ldots,n$ and $k=1,\ldots,m.$
\end{lemma}
\paragraph{Proof} See \cite{Baur1983TCo}. $\blacksquare$

An illustrative explanation of each method is followed by its formalization 
and the proof of numerical correctness. 

\subsection{Edge Splitting}

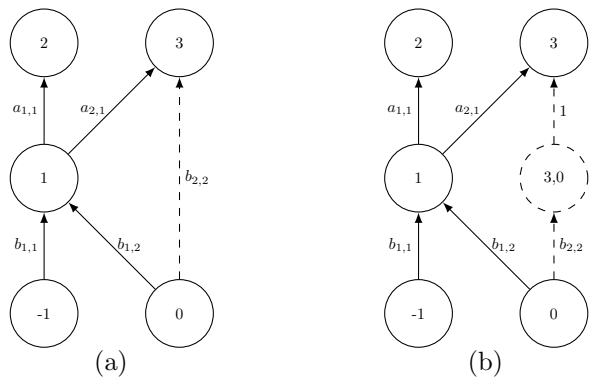
\begin{figure}
	\centering
\begin{tabular}{cc}
\begin{minipage}[c]{.35\textwidth}
	\centering
\begin{tikzpicture}[scale=.6, transform shape]
            \begin{pgfscope}
                \tikzstyle{every node}=[draw,circle,minimum size=1.5cm]
                    \node (00) at (0,0) {-1};
                    \node (01) at (3,0) {0};
                    \node (10) at (0,3) {1};
                    \node (20) at (0,6) {2};
                    \node (21) at (3,6) {3};
            \end{pgfscope}
            \begin{scope}[-latex]
                \draw (00) -- (10) node[midway,left] {$b_{1,1}$};
                    \draw[dashed] (01) -- (21) node[midway,right] {$b_{2,2}$};
                \draw (01) -- (10) node[midway,right] {$b_{1,2}$};
                \draw (10) -- (20) node[midway,left] {$a_{1,1}$};
                \draw (10) -- (21) node[midway,left] {$a_{2,1}$};
            \end{scope}
        \end{tikzpicture}
\end{minipage}
	&
\begin{minipage}[c]{.35\textwidth}
	\centering
\begin{tikzpicture}[scale=.6, transform shape]
            \begin{pgfscope}
                \tikzstyle{every node}=[draw,circle,minimum size=1.5cm]
                    \node (00) at (0,0) {-1};
                    \node (01) at (3,0) {0};
                    \node (10) at (0,3) {1};
                    \node[dashed] (11) at (3,3) {3,0};
                    \node (20) at (0,6) {2};
                    \node (21) at (3,6) {3};
            \end{pgfscope}
            \begin{scope}[-latex]
                \draw (00) -- (10) node[midway,left] {$b_{1,1}$};
                    \draw[dashed] (01) -- (11) node[midway,right] {$b_{2,2}$};
                \draw (01) -- (10) node[midway,right] {$b_{1,2}$};
                \draw (10) -- (20) node[midway,left] {$a_{1,1}$};
                    \draw[dashed] (11) -- (21) node[midway,right] {$1$};
                \draw (10) -- (21) node[midway,left] {$a_{2,1}$};
            \end{scope}
        \end{tikzpicture}
\end{minipage} \\
	(a) & (b)
\end{tabular}
\caption{Edge Splitting} \label{fig:es}
\end{figure}
Consider a minor modification of the sample program from Section~\ref{sec:1}
as
\begin{lstlisting}
for (int i=0;i<q;i++) {
  x[i%2]*=x[(i+1)%2];
  if (!(i%2)) x[(i+1)%2]=sin(x[(i+1)%2]);
}
\end{lstlisting}
For \lstinline{q=2}, the DAG in Figure~\ref{fig:es}(a) is recorded. 
The edge labeled with the
local partial derivative $b_{2,2}$ (=\lstinline{x[0]} for {i=1})
spans two layers making the DAG not layered. We aim to transform it into
a layered tripartite DAG representing the local Jacobian product $A \cdot B,$
where $A \equiv (a_{j,i})$ as well as $B \equiv (b_{j,i}).$ 
The $\partial_{j',i'}$ are replaced by $a_{j,i}$ or $b_{j,i}.$ 
This slight
modification in the notation is expected to make he upcoming examples easier
to follow.

Edge splitting makes the DAG layered
by inserting $l-1=j-i-1$ dummy vertices for all edges connecting vertices in 
layers $i$ and $j$
(here, $l=2-0=2$ yields one additional vertex labelled $3,0$ to mark it
as split vertex of edge $(0,3)$). One of the resulting new edges keeps the original label while
the others are labelled with ones as in Figure~\ref{fig:es}(b) ($l-1=1$ dummy vertex). The Jacobian
$F'$ turns out to be invariant under edge splitting as an immediate consequence of the 
chain rule of differentiation. For example,
$$
F'=
\begin{pmatrix}
        a_{1,1} & 0 \\
        a_{2,1} & 1
\end{pmatrix} 
\begin{pmatrix}
        b_{1,1} & b_{1,2} \\
        0 & b_{2,2}
\end{pmatrix}
=
\begin{pmatrix}
        a_{1,1} b_{1,1} & a_{1,1} b_{1,2} \\
        a_{2,1} b_{1,1} & a_{2,1} b_{1,2}+b_{2,2}
\end{pmatrix} \; .
$$
Invertibility of $F'$ implies invertibility of all elemental Jacobians in 
the resulting uniformly layered DAG.

\begin{definition}[Edge Splitting] \label{def:edgesplitting}
Let $G=(V,E)$ be a DAG as in Definition~\ref{def:dag}. 
	An {\em edge} $(i,j)$ with label $\lambda$ is {\em split} by replacing it with
two new edges $(i,k)$ and $(k,j),$ $k \not\in V,$ and labeling them
with $\partial_{k,i}=\lambda$ and $\partial_{j,k}=1$ (or vice versa), 
respectively.
\end{definition}
Consequences of the two choices for labelling the new edges are
the subject of ongoing investigations. Both alternatives are regarded
as equivalent for the purpose of the upcoming discussions.
\begin{lemma}
	Let a DAG $G$ induce an invertible Jacobian $F' \in \R^{n \times n}$ under 
the chain rule of differentiation as in Lemma~\ref{lem:cr}. 
\begin{enumerate}
\item The Jacobian is 
invariant under edge splitting as in Definition~\ref{def:edgesplitting}.
\item If repeated edge splitting yields a uniformly layered DAG as
in Definition~\ref{def:uldag}, then
the local Jacobians of each layer are invertible.
\end{enumerate}
\end{lemma}
\paragraph{Proof} 
\begin{enumerate}
\item Equation~(\ref{eqn:cr}) is invariant under edge splitting.
\item Let the uniformly layered DAG have depth $q.$ 
Equation~(\ref{eqn:cr}) yields $F'=F'_q \cdot \ldots \cdot F'_1$ and hence
		$$(F')^{-1}=(F'_q \cdot \ldots \cdot F'_1)^{-1}=(F'_1)^{-1} \cdot \ldots \cdot (F'_q)^{-1}$$ 
		implying
invertibility of all $F'_j,$ $j=1,\ldots,q.$
\end{enumerate} 
$\blacksquare$

Obviously, edge splitting terminates as soon as the DAG becomes layered.

\subsection{Preaccumulation}

\begin{figure}
	\centering
\begin{tabular}{cc}
\begin{tikzpicture}[scale=.6, transform shape]
            \begin{pgfscope}
                \tikzstyle{every node}=[draw,circle,minimum size=1.5cm]
                \node (n1) at (1.5,0) {$-2$};
                \node (n2) at (4.5,0) {$-1$};
                \node (n3) at (7.5,0) {$0$};
                \node (p1) at (0,3) {$1$};
                \node (p2) at (3,3) {$2$};
                \node (p3) at (6,3) {$3$};
                \node (p4) at (9,3) {$4$};
                \node (m1) at (1.5,6) {$5$};
                \node (m2) at (4.5,6) {$6$};
                \node (m3) at (7.5,6) {$7$};
            \end{pgfscope}
            \begin{scope}[-latex]
                \draw (n1) -- (p1) node[midway,left] {$b_{1,1}$};
                \draw (n1) -- (p3) node[midway,above] {};
                \draw (n2) -- (p2) node[near start,right] {$b_{2,2}$};
                \draw (n2) -- (p4) node[midway,above] {};
                \draw (n3) -- (p4) node[midway,above] {};
                    \draw[dashed] (p1) -- (m1) node[midway,left] {$a_{1,1}$};
                    \draw[dashed] (p2) -- (m1) node[midway,right] {$a_{1,2}$};
                \draw (p3) -- (m2) node[midway,left] {};
                \draw (p4) -- (m2) node[midway,above] {};
                \draw (p4) -- (m3) node[midway,above] {};
            \end{scope}
\end{tikzpicture} &
\begin{tikzpicture}[scale=.6, transform shape]
            \begin{pgfscope}
                \tikzstyle{every node}=[draw,circle,minimum size=1.5cm]
                \node (n1) at (0,0) {$-2$};
                \node (n2) at (3,0) {$-1$};
                \node (n3) at (6,0) {$0$};
                \node[dashed] (p12) at (0,3) {$1,2$};
                \node (p3) at (3,3) {$3$};
                \node (p4) at (6,3) {$4$};
                \node (m1) at (0,6) {$5$};
                \node (m2) at (3,6) {$6$};
                \node (m3) at (6,6) {$7$};
            \end{pgfscope}
            \begin{scope}[-latex]
                    \draw[dashed] (n1) -- (p12) node[midway,left] {$a_{1,1} b_{1,1}$};
                \draw (n1) -- (p3) node[midway,above] {};
                    \draw[dashed] (n2) -- (p12) node[near start,right] {$a_{1,2} b_{2,2}$};
                \draw (n2) -- (p4) node[midway,above] {};
                \draw (n3) -- (p4) node[midway,above] {};
                    \draw[dashed] (p12) -- (m1) node[midway,left] {$1$};
                \draw (p3) -- (m2) node[midway,above] {};
                \draw (p4) -- (m2) node[midway,above] {};
                \draw (p4) -- (m3) node[midway,above] {};
            \end{scope}
        \end{tikzpicture}
	\\
	(a) &(b)
\end{tabular}
\caption{Preaccumulation} \label{fig:preacc}
\end{figure}
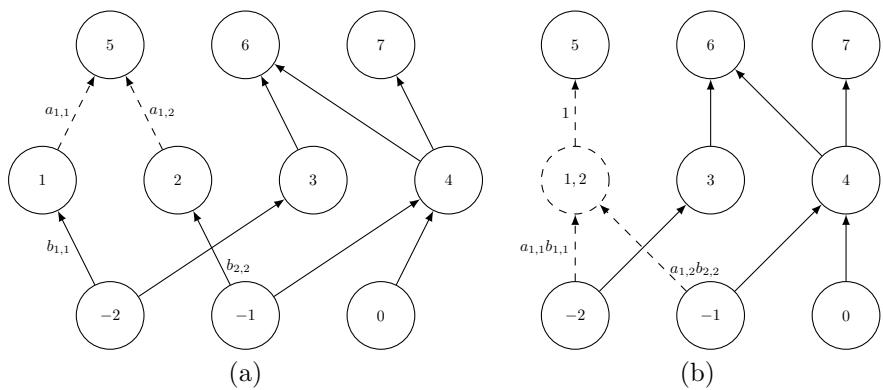
Making DAGs layered by edge splitting turned out to be rather straightforward.
Unfortunately, the resulting layered DAGs are rarely uniform. Moreover,
invertibility of the Jacobian implies that the sizes of all layers are
at least equal to $n.$
We propose 
{\em preaccumulation} of local sub-Jacobians followed by the 
decomposition of the resulting bipartite sub-DAGs into layered 
tripartite sub-DAGs as a method for making such DAGs uniformly layered.

Figure~\ref{fig:preacc}(a) shows a non-uniformly layered tripartite DAG $G$
with $n=m=3$ sources, respectively, sinks and $p=4$ intermediate 
vertices. 
Connected components of the bipartite sub-DAGs induce pure layered 
tripartite sub-DAGs
of $G$ as formally described in Definition~\ref{def:pltsdag}.
For example, the bipartite component spanned by the three vertices
$1,2,5$ induces a pure layered tripartite sub-DAG $\hat{G}$ of $G$ including
$-2$ and $-1.$ (Similarly, 
the bipartite component spanned by 
$-2,1,3$ yields 
the pure layered tripartite sub-DAG of $G$ including
$5$ and $6.$) Application of the chain rule of differentiation as in
Lemma~\ref{lem:cr} to $\hat{G}$ yields a local 
Jacobian $A \in \R^{m \times n}$ 
(here, $A \in \R^{1 \times 2}$ is the gradient of $v_5$ with 
respect to $v_{-2}$ and $v_{-1}$). 

Replacing $\hat{G}$ with the 
bipartite DAG of $A$ violates the requirement for $G$ 
to be layered. Decomposition of $A$ as
 $A=I_m \cdot A$ (or $A=A \cdot I_n$) 
takes care of this violation by substitution of a modified
pure layered tripartite sub-DAG $\tilde{G}$ as shown in Figure~\ref{fig:preacc}(b).
The number of intermediate vertices in $\tilde{G}$ is equal to 
$\tilde{p}=\min(n,m)$
implying a reduction of the number of intermediate vertices in $G$ if 
$\tilde{p}<p.$ Figure~\ref{fig:preacc} yields 
$1=\tilde{p}<p=2$ making the resulting DAG uniformly layered. 

Ultimately, preaccumulation always yields the 
uniformly layered DAG preserving all intermediate layers of size $n$ from
the original DAG.
In the worst case we get the bipartite DAG of $F'.$
Our objective is
to apply preaccumulation selectively as illustrated above. The resulting
uniformly layered DAG is made is meant to preserve the sparsity of all local 
Jacobians as well as possible. 

Decomposition of a bipartite graph (similarly, DAG) $G=(V,E)$ into its 
connected components turns out to be straightforward.
Each vertex belongs to exactly one connected component. Repeated (e.g., depth-first) searches starting from unassigned vertices yield a decomposition into 
connected components at a computational cost of ${\cal O}(|V|+|E|).$

The following definitions formalize the above followed by the proof of
the numerical correctness of preaccumulation in combination with bipartite DAG splitting.
\begin{definition}[Pure Layered Tripartite sub-DAG] \label{def:pltsdag}
A layered tripartite sub-DAG $G=(V,E),$ $V=(X,Z,Y),$ as in 
	Definition~\ref{def:ldag} is {\em pure} if $P(Z)=X$ and $S(Z)=Y.$
\end{definition}

\begin{definition}[Preaccumulation] \label{def:preaccumulation}
Let $G$ be a pure layered tripartite sub-DAG as in Definition~\ref{def:pltsdag} 
	representing a local Jacobian product $F'=F'_2 \cdot F'_1.$ {\em Preaccumulation} 
replaces $G$ by the bipartite sub-DAG representing $F'.$
\end{definition}
\begin{lemma}
	Let a DAG $G$ induce an invertible Jacobian $F' \in \R^{n \times n}$ under 
the chain rule of differentiation as in Lemma~\ref{lem:cr}. 
\begin{enumerate}
\item The Jacobian is 
invariant under preaccumulation as in Definition~\ref{def:preaccumulation}.
\item Preaccumulation terminates. 
\end{enumerate}
\end{lemma}
\paragraph{Proof} 
\begin{enumerate}
\item Equation~(\ref{eqn:cr}) is invariant under preaccumulation.
\item The cumulative sum of all paths connecting sources with sinks decreases
monotonically under preaccumulation. Hence, full preaccumulation yields the bipartite DAG that represents $F'.$
\end{enumerate} 
$\blacksquare$

\begin{definition}[Bipartite DAG Splitting] \label{def:bidagsplitting}
	A {\em bipartite DAG} $G$ representing a local
Jacobian $F'_j \in \R^{m \times n}$ 
	is {\em split} by replacing it with a layered tripartite DAG representing
the matrix products $F' \cdot I_n,$ if $n\leq m,$ 
or $I_m \cdot F',$ otherwise.
\end{definition}
The condition $n\leq m$ could be omitted if numerical correctness was our sole objective.
However, aiming for maximum sparsity we keep the number of newly generated edges
as low as possible. Further implications of the two choices are the subject 
of ongoing research.

\begin{lemma}
	Let a DAG $G$ induce an invertible Jacobian $F' \in \R^{n \times n}$ under 
the chain rule of differentiation as in Lemma~\ref{lem:cr}. 
\begin{enumerate}
\item The Jacobian is 
invariant under bipartite DAG splitting as in Definition~\ref{def:bidagsplitting}.
\item If repeated bipartite DAG splitting yields a uniformly layered DAG as
in Definition~\ref{def:uldag}, then
the local Jacobians of each layer are invertible.
\end{enumerate}
\end{lemma}
\paragraph{Proof} 
\begin{enumerate}
\item Equation~(\ref{eqn:cr}) is invariant under bipartite DAG splitting.
\item Let the uniformly layered DAG have depth $q.$ 
Equation~(\ref{eqn:cr}) yields $F'=F'_q \cdot \ldots \cdot F'_1$ and hence
$(F')^{-1}=(F'_1)^{-1} \cdot \ldots \cdot (F'_q)^{-1}$ implying
invertibility of all $F'_j,$ $j=1,\ldots,q.$
\end{enumerate} 
$\blacksquare$

\subsection{Combinatorics}

The combinatorial {\sc Matrix-Free Newton Step} problem to be formulated in its
decision version next turns out to be computationally intractable. 
\begin{problem}[\sc Matrix-Free Newton Step]
\item[]
        {\em \footnotesize INSTANCE:} An invertible DAG and an integer $K\geq 0.$ 
\item[]
        {\em \footnotesize QUESTION:} Can the corresponding Newton step be 
	evaluated with at most $K$ flops\footnote{floating-point operations}?
\end{problem}
\begin{theorem} \label{the}
	{\sc Matrix-Free Newton Step} is NP-complete.
\end{theorem}
\paragraph{Proof} The proof of NP-completeness of the 
{\sc Adjoint Computation} problem presented in \cite{Naumann2008OJA} reduces 
{\sc Ensemble Computation} \cite{Garey1979CaI} to uniformly layered
DAGs representing matrix chain products over diagonal Jacobians $\in \R^{n \times n}$ as
$$
F'=F'_q \cdot \ldots \cdot F'_1 \in \R^{n \times n}\; .
$$
A corresponding Newton step becomes equal to
\begin{align*}
(F')^{-1} \cdot y&= (F'_q \cdot \ldots \cdot F'_1)^{-1} \cdot y 
	= (F'_1)^{-1} \cdot \ldots \cdot (F'_q)^{-1} \cdot y \; .
\end{align*}
The ``accumulate first, then factorize'' method turns out to be superior as
$F'$ remains a diagonal matrix.
Its inversion adds a constant offset of $n$ 
(reciprocals) to the flop count of any given instance of {\sc Adjoint Computation}. A solution for {\sc Matrix-Free Newton Step} would hence solve 
{\sc Adjoint Computation} implying NP-hardness of the former. 
Moreover, a proposed solution is easily validated efficiently by counting the 
at most $(q+1)   n$ flops.
$\blacksquare$

The proof of NP-completeness of the {\sc Matrix-Free Newton Step} problem relies
entirely on algebraic dependences (equality in particular) amongst the elemental
partial derivatives. The structure of the underlying DAGs turns out to be 
trivial. On the other hand, the methods proposed for making DAGs uniformly 
layered are motivated by purely structural reasoning about the DAGs. All elemental 
partial derivatives were assumed to be mutually independent. At this stage, the
reduction used in the proof of Theorem~\ref{the} has no consequences on 
algorithms for the approximate solution of instances of {\sc Matrix-Free 
Newton Step}. In fact, similar statements apply to the related proof 
presented in \cite{Naumann2008OJA}.
Future research is expected to fill this gap.

\section{Conclusion} \label{sec:4}

The new matrix-free Newton method promises a significant reduction in the
computational cost of solving systems of nonlinear equations. This claim
is supported by run time measurements for problems with tridiagonal elemental
Jacobians. Further challenges need to be addressed for elemental Jacobians
with irregular sparsity patterns as well as for computationally expensive
residuals.

The matrix-free Newton method relies on the reversal of the DAG.
For very large problems the size of the DAG may exceed the available memory 
resources. 
The combinatorial {\sc DAG Reversal} problem asks for a distribution 
of the available storage such that the overall computational cost is
minimized. It is known to be NP-complete \cite{Naumann2008DRi}
as is the related {\sc Call Tree Reversal} problem \cite{Naumann2008CTR}.
Checkpointing methods \cite{Griewank1992ALG, Stumm2009MAf,Wang2009MRD} offer solutions for a 
variety of special cases. 

Factorization of elemental Jacobians with irregular sparsity patterns yields 
several well-known combinatorial problems in sparse linear algebra such as
{\sc Bandwidth} \cite{Papadimitriou1976TNo} or {\sc Directed Elimination Ordering} \cite{Rose1978Aao}.
Fill-in needs to be kept low for our method to outperform other
state-of-the-art solutions. A rich set of results from past and ongoing efforts
in this highly active area of research can be built on.

Another class of potential targets are surrogates with banded layers for 
computationally expensive nonlinear residuals obtained by machine learning 
(ML). 
While the training of such models can be challenging, the inversion of their
Jacobian is guaranteed to benefit from the new method. Ongoing research aims
to explore a potential extension to ML models with layers exhibiting other suitable sparsity patterns, e.g, with triangular elemental Jacobians.

\subsection*{Acknowledgement}

The experiments reported on in Figure~\ref{fig:3} were performed by
Gero Kauerauf as part of his M.Sc. thesis project at RWTH Aachen University. 

\bibliographystyle{siamplain}
\bibliography{paper}

\end{document}